\documentclass[12pt]{amsart}
\usepackage{graphicx}
\usepackage{amssymb}
\usepackage{amsfonts}
\usepackage{amsmath,amscd}
\usepackage{array}
\usepackage{rotating}
\usepackage{epsfig}
\usepackage[matrix,frame]{xypic}

\newtheorem{example}{Example}[section]

\newtheorem{theorem}[example]{Theorem}

\newtheorem{corollary}[example]{Corollary}

\newtheorem{proposition}[example]{Proposition}

\newtheorem{lemma}[example]{Lemma}

\def\Proof{\noindent \it Proof -- \rm}
\def\qed{\hspace{3.5mm} \hfill \vbox{\hrule height 3pt depth 2 pt width 2mm}
\bigskip}

\def\QSym{{\it QSym}}          % QSym
\def\FQSym{{\bf FQSym}}        % permutations
\def\PBT{{\bf PBT}}            % arbres binaires planaires
\def\BT{{\bf BT}}
\def\CBT{{\bf CBT}}
\def\WQSym{{\bf WQSym}}        % Mots initiaux 
\def\maj{{\rm imaj\,}}
\def\maj{{\rm maj\,}}

\def\P{{\bf P}}
\def\M{{\bf M}}
\def\pack{{\rm pack\,}}
\def\PW{{\rm PW}}
\def\Std{{\rm Std}}     % standardisation
\def\<{\langle}
\def\>{\rangle}
\def\EE{{\mathbb E}\, } % alphabet K

\def\F{{\bf F}}         % F de FQSym
\def\G{{\bf G}}         % G de FQSym^*
\def\SG{{\mathfrak S}}  % groupe symetrique
    
\def\B{{\bf B}}

\def\K{{\mathbb K}}
\def\X{{\bf X}}
\def\B{{\bf B}}
\def\T{{\mathcal T}}

\def\maj{{\rm maj\,}}
\def\imaj{{\rm imaj\,}}

\def\shuff#1#2{\mathbin{
\hbox{\vbox{ \hbox{\vrule \hskip#2 \vrule height#1 width 0pt
}%
\hrule}%
\vbox{ \hbox{\vrule \hskip#2 \vrule height#1 width 0pt
\vrule }%
\hrule}%
}}}

\long\def\psboxit#1#2{%
\begingroup\setbox0=\hbox{#2}%
\dimen0=\ht0 \advance\dimen0 by \dp0%
    \hbox{%
    \copy0%
    }%hbox
\endgroup%
}%psboxit
\def\SetTableau#1#2#3#4{%
  \gdef\Tabvrule{\vrule\vrule width-0.4pt}
  \gdef\Tabhrule{\hrule\hrule height-0.4pt}  
  \gdef\Tabstrut{\vrule height#1 depth#2 width0pt\relax}
  \gdef\Tabbox##1{\hbox to #3{\hskip0.4pt\hfill\Tabstrut$#4##1$\hfill}}
} %setTableau
\def\Case#1{\vcenter{\Tabhrule%
                   \hbox{\Tabvrule\Tabbox{#1}\Tabvrule}\Tabhrule}}
\def\GenTab#1{\vcenter{\halign{&$\Case{##}$\cr#1}}\egroup}
\def\Tableau{%         
  \bgroup%
  \let\ =\omit%
  \let\\=\cr%
  \offinterlineskip\GenTab}

\def\qbin#1#2{\begin{bmatrix} #1 \\ #2\end{bmatrix}}

\def\shuf{{\mathchoice{\shuff{7pt}{3.5pt}}%
{\shuff{6pt}{3pt}}%
{\shuff{4pt}{2pt}}%
{\shuff{3pt}{1.5pt}}}}%
\def\shuffle{\,\shuf\,}

\title[Trees and combinatorial Hopf algebras]%
{Trees, functional equations, and combinatorial Hopf algebras}

\author[F. Hivert, J.-C.~Novelli, and J.-Y.~Thibon]%
{Florent Hivert, Jean-Christophe Novelli, and Jean-Yves Thibon}

\address[Hivert]{LITIS, Universit\'e de Rouen ; Avenue de l'universit\'e ;
76801 Saint \'Etienne du Rouvray, France\\}

\address[Novelli and Thibon]{Institut Gaspard Monge, Universit\'e de
Marne-la-Vall\'ee \\
5, Boulevard Descartes \\Champs-sur-Marne \\77454 Marne-la-Vall\'ee cedex 2 \\
FRANCE}
\email[Florent Hivert]{hivert@univ-rouen.fr}
\email[Jean-Christophe Novelli]{novelli@univ-mlv.fr}
\email[Jean-Yves Thibon]{jyt@univ-mlv.fr}

\date{}

\begin{document}

\begin{abstract}
One of the main virtues of trees is to represent formal solutions of various
functional equations which can be cast in the form of fixed point problems.
Basic examples include differential equations and functional (Lagrange)
inversion in power series rings. When analyzed in terms of combinatorial Hopf
algebras, the simplest examples yield interesting algebraic identities or
enumerative results.
\end{abstract}

\maketitle

%{\footnotesize
%\tableofcontents
%}

%%%%%%%%%%%%%%%%%%%%%%%%%%%%%%%%%%%%%%%%%%%%%%%%%%%%%%%%%%%%%%%%%%%%%%%%%%%%%%%
%%%%%%%%%%%%%%%%%%%%%%%%%%%%%%%%%%%%%%%%%%%%%%%%%%%%%%%%%%%%%%%%%%%%%%%%%%%%%%%
%%%%%%%%%%%%%%%%%%%%%%%%%%%%%%%%%%%%%%%%%%%%%%%%%%%%%%%%%%%%%%%%%%%%%%%%%%%%%%%
\section{Introduction}

Let $R$ be an associative algebra, and consider the functional
equation for the power series $x\in R[[t]]$ 
\begin{equation}\label{eqbin}
x = a + B(x,x)
\end{equation}
where $a\in R$ and $B(x,y)$ is a bilinear map with values in
$R[[t]]$, such that the valuation of $B(x,y)$ is strictly greater
than the sum of the valuations of $x$ and $y$. Then, (\ref{eqbin})
has a unique solution
\begin{equation}\label{solbin}
x = a + B(a,a) + B(B(a,a),a)+ B(a,B(a,a))+ \cdots =  \sum_{T\in\CBT} B_T(a)
\end{equation}
where $\CBT$ is the set of (complete) binary trees, and for a tree $T$,
$B_T(a)$ is the result of evaluating the expression formed by labeling
by $a$ the leaves of $T$ and by $B$ its internal nodes.

Of course, the same can be done with $m$-ary trees, or more generally
with {\em plane trees}. 
We are in particular interested in those counted by the little Schr\"oder numbers,
that is, plane trees without vertex of arity $2$ \cite[A001003]{Slo}, 
which solve equations of the
form
\begin{equation}\label{eqplane}
x=a+\sum_{n\ge 2}F_n(x,x,\ldots,x)
\end{equation}
each $F_n$ being an $n$-linear operation.

All this is well-known and rather trivial.
However, the simplest example has still something
to tell us.
Consider the differential equation (for $x\in\K[[t]]$)
\begin{equation}\label{eqbindif}
\frac{dx}{dt} = x^2\,,\ \ x(0)=1\,.
\end{equation}
Its solution is obviously $x=(1-t)^{-1}$, but let us ignore this for
the moment, and recast it as a fixed point problem
\begin{equation}\label{eqbinint}
x = 1 + \int_0^t x^2(s) ds = 1 + B(x,x)\,,
\end{equation}
where $B(x,y)=\int_0^t x(s)y(s)ds$. Then, for a binary tree $T$ with $n+1$
leaves, $B_T(1)$ is the monomial obtained by putting $1$ on each leaf
and integrating at each internal node the product of the evaluations of
its subtrees:
\entrymodifiers={+<4pt>}
\begin{equation}
\vcenter{\xymatrix@C=2mm@R=2mm{
*{}   &  *{} & *{} & {t^4/8}\ar@{-}[drr]\ar@{-}[dll] \\
*{}   & {t}\ar@{-}[dr]\ar@{-}[dl]  & *{} & *{}   & *{} &
{t^2/2}\ar@{-}[dl]\ar@{-}[ddr] \\
{1}   & *{}  & {1} & *{}   & {t}\ar@{-}[dl]\ar@{-}[dr]\\
*{}   & *{}  & *{} & {1}   & *{} & {1} & {1}\\
      }}
\end{equation}

One can observe that 
\begin{equation}
B_T(1) = c_{T'} \frac{t^n}{n!}\,,
\end{equation}
where $T'$ is the incomplete binary tree with $n$ nodes obtained by
removing the leaves of $T$, and $c_{T'}$ is the number of permutations
$\sigma\in\SG_n$ whose decreasing tree has shape $T'$.
Indeed, $c_{T'}$ is explicitly given by a {\em hook length formula}
\cite{Kn}, which can be compared with the easily obtained closed form
for $B_T(1)$. The {\em hook lengths} of $T'$ are the number of nodes
of all the subtrees
\entrymodifiers={+<4pt>[o][F]}
\begin{equation}
    \vcenter{\xymatrix@C=2mm@R=2mm{
        *{}  & {4}\ar@{-}[dr]\ar@{-}[dl] \\
        {1}  & *{} & {2}\ar@{-}[dl] \\
        *{}  & {1}\\
      }}
\end{equation}
and $c_{T'}$ is $n!$ over the product of the hook lengths,
here $4!/(4\cdot 2\cdot 1\cdot 1)=3$, the corresponding decreasing trees
being
\entrymodifiers={+<4pt>}
\begin{equation}
    \vcenter{\xymatrix@C=2mm@R=2mm{
        *{}  & {4}\ar@{-}[dr]\ar@{-}[dl] \\
        {1}  & *{} & {3}\ar@{-}[dl] \\
        *{}  & {2}\\
      }}
\qquad
    \vcenter{\xymatrix@C=2mm@R=2mm{
        *{}  & {4}\ar@{-}[dr]\ar@{-}[dl] \\
        {2}  & *{} & {3}\ar@{-}[dl] \\
        *{}  & {1}\\
      }}
\qquad
    \vcenter{\xymatrix@C=2mm@R=2mm{
        *{}  & {4}\ar@{-}[dr]\ar@{-}[dl] \\
        {3}  & *{} & {2}\ar@{-}[dl] \\
        *{}  & {1}\\
      }}
\end{equation}

Our starting point will be the following question: can one use this observation
to derive the hook length formula for binary trees, and if yes, can we
use the same method to obtain more interesting results ?

For this, we have to lift our problem to the
combinatorial Hopf algebra of Free quasi-symmetric functions $\FQSym$. 
We can then derive in the same way the $q$-hook length formulas of Bj\"orner
and Wachs
\cite{BW1,BW2}. The case of plane trees can be dealt with in the same way,
the relevant Hopf algebra being there $\WQSym$, the Word Quasi-Symmetric
invariants (or quasi-symmetric functions in noncommutative variables), and
here the resulting formula is believed to be new. 
Finally, we give new proofs  of some identities
of Postnikov \cite{apost} and Du-Liu \cite{DuLiu} by relating these
to appropriate functional equations.

{\small
\subsection*{Notations} The symmetric group is denoted by $\SG_n$.
The standardized $\Std(w)$ of a word $w$ of length $n$ is
the permutation obtained by iteratively scanning $w$ from
left to right, and labelling $1,2,\ldots$ the occurrences of its
smallest letter, then numbering the occurrences of the next one, and
so on. All algebras are over a field $\K$ of characteristic 0.
}

\section{Free quasi-symmetric functions and hook length formulas}

\subsection{A derivation of $\FQSym$}

Recall from \cite{NCSF6} that for an (infinite) totally ordered alphabet $A$,
$\FQSym(A)$ is the subalgebra of $\K\<A\>$ spanned by the polynomials
\begin{equation}
\G_\sigma(A)=\sum_{\Std(w)=\sigma}w
\end{equation}
the sum of all words in $A^n$ whose standardization is the permutation
$\sigma\in\SG_n$. The multiplication rule is, for $\alpha\in\SG_k$ and
$\beta\in\SG_l$,
\begin{equation}\label{multG}
\G_\alpha \G_\beta=\sum_{\gamma\in\SG_{k+l};\,\gamma=u\cdot v\atop
\Std(u)=\alpha,\Std(v)=\beta}\G_\gamma\,.
\end{equation}
This sum has ${k+l\choose k}$
terms. Hence, the linear map
\begin{equation}
\phi:\ G_\sigma \longmapsto \frac{t^n}{n!} \quad (\sigma\in\SG_n)
\end{equation}
is a homomorphism of algebras $\FQSym\rightarrow \K[[t]]$.
It is convenient to introduce the notation
$\F_\sigma=\G_{\sigma^{-1}}$ and a scalar product
satisfying $\<\F_\sigma,\G_\tau\>=\delta_{\sigma,\tau}$.
As a graded bialgebra, $\FQSym$ is self-dual, and its coproduct
$\Delta$ satisfies $\<FG,H\>=\<F\otimes G,\Delta H\>$.

Let $\partial$ be the linear map defined by
\begin{equation}
\partial\G_\sigma=\G_{\sigma'}
\end{equation}
where $\sigma'$ is the permutation whose word is obtained by erasing 
the letter $n$ in $\sigma\in\SG_n$. Obviously,
\begin{equation}
\phi(\partial F) = \frac{d}{dt}\phi(F)
\end{equation}
for all $F\in\FQSym$, and moreover:

\begin{proposition}\label{partial}
The map $\partial$ is a derivation of $\FQSym$. It is the
adjoint of the linear map $F\mapsto F\cdot \F_1$.
\end{proposition}

\Proof By definition, $\<\partial\G_\sigma,\F_\tau\>=\delta_{\sigma',\tau}$
is equal to 1 if $\sigma$ occurs in $\tau\shuffle n$ and to $0$ otherwise.
Hence, 
\begin{equation}
\<\partial\G_\sigma,\F_\tau\>=\<\G_\sigma,\F_\tau\F_1\>\,
\end{equation}
whence the second part of the proposition. Now, $\F_1$
is a primitive element, so that $\partial$ is a derivation.
\qed

%%%%%%%%%%%%
The Leibniz relation
\begin{equation}
\partial (\G_\alpha \G_\beta)=\partial\G_\alpha\cdot\G_\beta+\G\alpha\cdot\partial\G_\beta
\end{equation}
can be interpreted in terms of the {\em dendriform structure} of $\FQSym$.
Recall \cite{LR} that the product $\G_\alpha\G_\beta$ can be split
into two parts (the dendriform operations)
\begin{equation}
\G_\alpha\G_\beta =\G_\alpha\prec\G_\beta+\G_\alpha\succ\G_\beta\,,
\end{equation}
\begin{equation}
\G_\alpha\prec\G_\beta=\sum_{\gamma=u\cdot v,\ \max(u)>\max(v)\atop
\Std(u)=\alpha,\Std(v)=\beta}\G_\gamma\,,
\end{equation}
\begin{equation}
\G_\alpha\succ\G_\beta=
\sum_{\gamma=u\cdot v,\ \max(u)\le \max(v)\atop
\Std(u)=\alpha,\Std(v)=\beta}\G_\gamma\,.
\end{equation}
Then,
\begin{equation}
\partial (\G_\alpha\prec \G_\beta)=\partial\G_\alpha\prec\G_\beta\,,\quad
\partial (\G_\alpha\succ \G_\beta)=\G_\alpha\succ\partial\G_\beta\,.
\end{equation}
It will be convenient to consider the half products as also defined
on permutations, so that their sum is then the convolution $\alpha *\beta$.

%%%%%%%%%%%%%
\subsection{A differential equation in $\FQSym$}

It follows from Proposition \ref{partial} that
if we set $\X=(1-\G_1)^{-1}$, we have
\begin{equation}
\partial \X = \X^2
\end{equation}
with $\X_0=1$ (constant term), and $\phi(\X)=(1-t)^{-1}$.

Note that thanks to the multiplication formula (\ref{multG}),
\begin{equation}
\X = \sum_\sigma \G_\sigma = \sum_{w\in A^*} w
\end{equation}
is the sum of all permutations (interpreted as $\G$'s),
that is, the sum of all words. If we can lift to $\FQSym$
the scalar bilinar map $B(x,y)=\int_0^tx(s)y(s)ds$, it
will also be interpretable as the sum of all complete binary trees.

\subsection{The bilinear map} The required map is given by
a simple operation, already introduced in \cite{NCSF6}, precisely
with the aim of providing a better understanding of the Loday-Ronco
algebra \cite{LR} of planar binary trees.

For $\alpha\in\SG_k$, $\beta\in\SG_l$, and $n=k+l$, set
\begin{equation}
\B(\G_\alpha,\G_\beta)=\sum_{\gamma=u(n+1)v\atop
\Std(u)=\alpha,\Std(v)=\beta}\G_\gamma\,.
\end{equation}
Clearly,
\begin{equation}
\partial \B(\G_\alpha,\G_\beta)=\G_\alpha\G_\beta\,,
\end{equation}
and our differential equation is now equivalent to the fixed
point problem
\begin{equation}\label{eqbinX}
\X = 1 + \B(\X,\X)\,.
\end{equation}

\begin{theorem}
In the binary tree solution~{\rm (\ref{solbin})} of {\rm (\ref{eqbinX})}, 
\begin{equation}\label{BT1}
\B_T(1)=\sum_{\T(\sigma)=T}\G_\sigma,
\end{equation}
where $\T(\sigma)$ denotes the shape of the decreasing tree of
the permutation $\sigma$. In particular, $\B_T(1)$ 
coincides with $\P_T$, the natural basis of the Loday-Ronco algebra 
(in the notation of \cite{HNT}).
\end{theorem}

\Proof By induction on the number $n$ of internal nodes of $T$.
For $n=1$ the result is obvious, and if $n>1$, 
$$ \B_T(1)=\B(\B_{T'}(1),\B_{T''}(1))\,,$$
where $T'$ and $T''$ are the left an right subtrees of $T$.
Hence, $\B_T(1)$ is the sum of the $\G_\sigma$ for $\sigma=\alpha n\beta$
such that $\G_{\Std(\alpha)}$ occurs in $\B_{T'}(1)$ and $\G_{\Std(\beta)}$
occurs in $\B_{T''}(1)$. Since we have assummed that (\ref{BT1}) holds
for $T'$ and $T''$, this implies that it holds for $T$ as well.
\qed

\begin{corollary}[The hook length formula]
The number of permutations whose decreasing tree has shape $T$ is
\begin{equation}
\frac{n!}{\prod_{v\in T} h_v}\,,
\end{equation}
where for a vertex $v$ of $T$, $h_v$ is the number
of nodes of the subtree with root $v$.
\end{corollary}

\subsection{The $q$-hook length formula}

Recall that under the $q$-specialization
\begin{equation}
A= \frac{1}{1-q} := \{\ldots <q^n<q^{n-1}<\ldots <q<1\}
\end{equation}
we have \cite[(125)]{NCSF2}
\begin{equation}
\G_\sigma\left(\frac{1}{1-q}\right) = \frac{q^{\imaj(\sigma)}}{(q)_n}
\end{equation}
where $\imaj(\sigma)=\maj(\sigma^{-1})$, 
$\maj(\sigma)$ is the classical major index
(sum of the descents) of
$\sigma\in\SG_n$ and
$(q)_n=(1-q)(1-q^2)\cdots (1-q^n)$.

Hence, the map
\begin{equation}
\phi_q(\G_\sigma)=
\frac{q^{\imaj(\sigma)}t^n}{[n]_q!}=(t(1-q))^n\G_\sigma\left(\frac{1}{1-q}\right)
\end{equation}
is a homomorphism of algebras.
The image of (\ref{eqbinX}) under $\phi_q$ reads 
\begin{equation}
x = 1 + B_q(x,x)
\end{equation}
where the bilinear map is now a $q$-integral
\begin{equation}
B_q(f,g)=\int_0^t d_qs f(s)g(qs)\,,
\end{equation}
where the $q$-integral is defined by
\begin{equation}
\int_0^t s^n d_qs = \frac{t^{n+1}}{[n+1]_q}.
\end{equation}
To show this, we have to compute $\phi_q(\B(\G_\alpha,\G_\beta))$.
\begin{lemma} Let $\alpha\in\SG_k$, $\beta\in\SG_l$.
The inverse major index  is distributed over the half-products according to
\begin{equation}
\sum_{\gamma\in\alpha\succ\beta}q^{\imaj(\gamma)}=
q^{\imaj(\alpha)+\imaj(\beta)}\qbin{k+l-1}{l-1}_q\,,
\end{equation}
and
\begin{equation}
\sum_{\gamma\in\alpha\prec\beta}q^{\imaj(\gamma)}=
q^{\imaj(\alpha)+\imaj(\beta)+l}\qbin{k+l-1}{l}_q\,.
\end{equation}
\end{lemma}
 
\Proof Straightformward by induction on $n=k+l$. \qed

From this, on deduces immediately
\begin{equation}
\sum_{\gamma=u\cdot (n+1)\cdot v\atop\Std(u)=\alpha,\Std(v)=\beta}
q^{\imaj(\gamma)}=
q^{\imaj(\alpha)+\imaj(\beta)+l}\qbin{k+l}{k}_q\,,
\end{equation}
which in turn implies the following: 
%%%%%%%%%%%%ICI%%%%%%%%%%%%%%%%%%%%%%
\begin{lemma}
If $f(t)=\phi_q(F)$ and $g(t)=\phi_q(G)$, then
\begin{equation}
\phi_q(\B(F,G))=\int_0^t d_qs f(s)g(qs)\,.
\end{equation}
\end{lemma}

\begin{corollary}[The $q$-hook length formula of \cite{BW1}]
\label{qhook}
The inverse major index polynomial of the set of permutations
whose decreasing tree has shape $T$ is
\begin{equation}
\sum_{{\mathcal T}(\sigma)=T}q^{\imaj(\sigma)}
=
[n]_q! \prod_{v\in T}\frac{q^{\delta_v}}{[h_v]_q}\,,
\end{equation}
where $\delta_v$ is the number of nodes in the right
subtree of $v$.
 
\end{corollary}

\subsection{Another approach}
It has been observed in \cite{NCSF6} that $\FQSym$ had
a natural $q$-deformation, obtained by replacing the ordinary
shuffle $\shuffle$ by the $q$-shuffle $\shuffle_q$ in the product
formula for the basis $\F_\sigma$.
That is, $\FQSym_q$ is the algebra with basis $\F_\sigma= \G_{\sigma^{-1}}$ and
product rule
\begin{equation}
\F_\alpha \F_\beta = \sum_\gamma (\gamma|\alpha\shuffle_q\beta[k])\F_\gamma
=\sum_\gamma (\gamma|\alpha\shuffle\beta[k]) q^{l(\gamma)-l(\beta)-l(\alpha)}\F_\gamma
\end{equation}
where $(\gamma|f)$ means the coefficient of $\gamma$ in $f$, $k$
is the length of $\alpha$ and $\beta[k]=(\beta_1+k)\cdots(\beta_l+k)$, (the
shifted word), $l(\sigma)$ being the number of inversions of $\sigma$.

Then, the map $\phi_q:\ \FQSym_q\rightarrow \K[[t]]$ defined by
\begin{equation}
\phi_q(\G_\sigma) = \frac{t^n}{[n]_q!}
\end{equation}
is a homomorphism of algebras.

One has now
\begin{equation}
\phi_q(\partial F)= D_q\phi_q(F)
\end{equation}
where $D_q$ is the $q$-derivative
\begin{equation}
D_qf(t)=\frac{f(qt)-f(t)}{qt-t}\,.
\end{equation}
In $\FQSym_q$, $\partial$ is not anymore a derivation, but satisfies
\begin{equation}
\partial(FG)=\partial F(A)\cdot  G(qA)+F(A)\cdot \partial G(A) 
\end{equation}
so that the noncommutative functional equation is now
\begin{equation}
\partial \X(A)= \X(A)\X(qA)\,,\ \X_0=1
\end{equation}
and its one-variable projection under $\phi_q$ is
\begin{equation}
D_q x(t)=x(t)x(qt)\,,\ x(0)=1\,.
\end{equation}
This is equivalent to
\begin{equation}
x = 1+B_q(x,x)
\end{equation}
where we have again
\begin{equation}
B_q(x,y) = \int_0^td_qs\, x(s)y(qs)\,.
\end{equation}

\begin{theorem}[$q$-hook length formula for inversions \cite{BW2}]
The inversion polynomial of the set of permutations having
a decreasing tree of shape $T$  is given by the
same hook length formula as for the inverse major index,
\begin{equation}
\sum_{{\mathcal T}(\sigma)=T}q^{l(\sigma)}
=
[n]_q! \prod_{v\in T}\frac{q^{\delta_v}}{[h_v]_q}\,,
\end{equation}
In particular $\imaj$ and $l$ are equidistributed on these
sets.
\end{theorem}
This is a refinement of a classical result of Foata and Sch\"utzenberger.

\section{Word quasi-symmetric functions and plane trees}

To interpret (\ref{eqplane}), we need to work in $\WQSym$,
the algebra of Word Quasi-Symmetric functions,
which contains an algebra of plane trees (the free dendriform
trialgebra on one generator \cite{LRtri}) in the same way as
$\FQSym$ contains an algebra of binary trees \cite{NT-cras}.

The basis elements $\M_u$ of $\WQSym$ are labeled by
packed words $u$, or if one prefers, surjections
$[n]\rightarrow [k]$, set compositions, 
or facets of the permutohedron \cite{Ch}.
These objects are counted by the ordered Bell numbers
\cite[A000262]{Slo}.
There is a canonical way to associate a plane tree
to such an object \cite{NT-cras}, and the sums over the
fibers of this map span a Hopf subalgebra of $\WQSym$.
Hence, we need to define on $\WQSym$ an analogue
of our derivation $\partial$ of $\FQSym$.

Recall that a word $w$ over the aphabet of positive
integers is said to be {\em packed} if the set of letters
occuring in $w$ is an initial interval $[a_1,a_k]$ of the
alphabet $A$. 
The {\em packed word} $u=\pack(w)$ associated to a word $w\in A^*$ is
obtained by the following process. If
$b_1<b_2<\ldots <b_k$ are the letters occuring in $w$, $u$ is the
image of $w$ by the semigroup homomorphism $b_i\mapsto a_i$.
For example, $\pack(34364)=12132$.
A word $u$ is said to be {\em packed} if $\pack(u)=u$.
To such a word is associated a polynomial $\M_u$, 
defined as the sum of all words $w$ such that $\pack(w)=u$.

The product on $\WQSym$ is given by
\begin{equation}
\label{prodG-wq}
\M_{u'} \M_{u''} = \sum_{u \in u'\star u''} \M_u\,,
\end{equation}
where the \emph{convolution} $u'\star u''$ of two packed words
is defined as
\begin{equation}
u'\star u'' = \sum_{v,w ;
u=v\cdot w\,\in\,\PW, \pack(v)=u', \pack(w)=u''} u\,.
\end{equation}
For example,
\begin{equation}
\M_{11} \M_{21} =
\M_{1121} + \M_{1132} + \M_{2221} + \M_{2231} + \M_{3321}.
\end{equation}

The coproduct can  be defined  by the usual
trick of noncommutative symmetric functions, considering the
alphabet $A$ as an ordered sum of two mutually commuting alphabets
$A'\hat+A''$. First, by direct inspection, one finds that 
\begin{equation}\label{coprodM}
\M_u(A'\hat+A'') = \sum_{0\leq k\leq \max(u)}
\M_{(u|_{[1,k]})}(A') \M_{\pack(u|_{[k+1,\max(u)})}(A''),
\end{equation}
where $u|_{B}$ denote the subword obtained by restricting $u$ to the
subset $B$ of the alphabet.

For a packed word $u$, let $u'$ be the word obtained from
$u$ by erasing all the occurences of the maximal letter
$m=\max(u)$, e.g., $(5211354)'=21134$.
Now, define
a linear map $\delta$ by
\begin{equation}\label{delta}
\delta\M_u=\M_{u'}\,.
\end{equation}
This is not anymore a derivation, but rather a finite
difference operator: indeed, it follows from (\ref{coprodM})
that
\begin{equation}
\delta\M_u(A)=\M_u(A\hat+1)-\M_u(A)\,,
\end{equation}
where $A\hat+1$ is the ordered sum of $A$ and $\{1\}$ (the
scalar 1, so that $\M_u(1)=1$ if $u$ is of the form
$11\cdots 1$, and is 0 otherwise).
Alternatively, $\delta$ is the adjoint of the right multiplication by
$\sum_{n\geq1} \M^*_{1^n}$, where $\M_u^*$ is the dual basis of $\M_u$.

This implies that $\delta$ satisfies
\begin{equation}\label{fdo}
\delta(FG)=(\delta F)G+(\delta F)(\delta G)+F(\delta G)\,,
\end{equation}
but this formula can be refined in terms of the tridendriform
structure of $\WQSym$ \cite{NT-cras}.
Indeed, it is known that
$\WQSym^+$ is a sub-dendriform trialgebra of $\K\<A\>^+$,
the partial products being given by
\begin{equation}
\M_{w'} \prec \M_{w''} =
\sum_{w=u\cdot v\in w'\star w'', |u|=|w'| ; \max(v)<\max(u)}
\M_w,
\end{equation}
\begin{equation}
\M_{w'} \circ \M_{w''} =
\sum_{w=u\cdot v\in w'\star w'', |u|=|w'| ; \max(v)=\max(u)}
\M_w,
\end{equation}
\begin{equation}
\M_{w'} \succ \M_{w''} =
\sum_{w=u\cdot v\in w'\star w'', |u|=|w'| ; \max(v)>\max(u)}
\M_w\,r,.
\end{equation}
and it follows from the multiplication rule
(\ref{prodG-wq}) that
\begin{equation}
\delta(F\prec G)=(\delta F)G\,,\ \delta(F\circ G) =(\delta F)(\delta G)\,,\
\delta(F\succ G) = F(\delta G)\,.
\end{equation}

Now, let
\begin{equation}
\X=(1-q\M_1)^{-1}=\sum_u q^{|u|}\M_u =\sum_w q^{|w|}w\,.
\end{equation}
It follows from (\ref{fdo}) that
\begin{equation}\label{eqdw}
\delta\X=q\X^2(1-q\X)^{-1}=\sum_{n\ge 2}q^{n-1}\X^n\,.
\end{equation}

For packed words $u_1,\ldots,u_k$, define 
\begin{equation}\label{Fk}
\F_k(\M_{u_1},\ldots,\M_{u_k})=
\sum \M_w
\end{equation}
where the sums runs over packed words $w$ such that
\begin{equation}\label{wFk}
w=w_1 m w_2 m \cdots m w_k\,,\quad \pack(w_i)=u_i\,,\
m=\max(w_1,\ldots,w_k)+1\,.
\end{equation}
For example,
\begin{equation}
\F_2(\M_{11}, \M_{21}) =
\M_{11321} + \M_{11432} + \M_{22321} + \M_{22431} + \M_{33421}.
\end{equation}
Then, obviously,
\begin{equation}\label{eqfw}
\X=1+\sum_{n\ge 2}q^{n-1}\F_n(\X\ldots,\X)\,.
\end{equation}
which does indeed give back (\ref{eqdw}), since
\begin{equation}
\delta \F_k(\M_{u_1},\ldots,\M_{u_k})=\M_{u_1}\cdots\M_{u_k}\,.
\end{equation}

It follows from (\ref{prodG-wq}) that the linear map
$\psi:\ \WQSym\rightarrow\K[[t]]$ defined by
\begin{equation}
\psi(\M_u)={t\choose\max(u)}
\end{equation}
is a homomorphism of algebras. Moreover, it maps $\delta$
over the finite difference operator
\begin{equation}
\psi(\delta F)=\Delta\psi(F)
\end{equation}
where $\Delta f(t)=f(t+1)-f(t)$.
Hence, the images of (\ref{eqdw}) and (\ref{eqfw}) by $\psi$
are
\begin{eqnarray}
\Delta x=\sum_{n\ge 2}q^{n-1} x^n\,\\
x=1+\sum_{n\ge 2}q^{n-1}F_n(x,x,\ldots,x)\,,
\end{eqnarray}
where
\begin{equation}
F_n(x_1,\ldots,x_n)=\Sigma_0^t x_1(s)x_2(s)\cdots x_n(s)\delta s,
\end{equation}
the discrete integral being defined by
\begin{equation}
\Sigma_0^t f(s)\delta s =\sum_{i=0}^{t-1}f(i)\,.
\end{equation}
%%%%%%%%

\def\TT{{\mathcal T}} % Arbre plan associe a un mot

The realization of the free dendriform trialgebra
given in \cite{NT-cras} involves the following construction.
With any word $w$ of length $n$, associate a plane tree $\TT(w)$ with $n+1$
leaves, as follows: if $m=\max(w)$ and if $w$ has exactly $k-1$ occurences of
$m$, write
\begin{equation}\label{wT}
w=v_1\,m\,v_2\cdots v_{k-1}\,m\,v_k\,,
\end{equation}
where the $v_i$ may be empty. Then, $\TT(w)$ is the tree obtained by
grafting the subtrees $\TT(v_1),\TT(v_2),\ldots,\TT(v_k)$ (in this order)
on a common root, with the initial condition $\TT(\epsilon)=\emptyset$
for the empty word.
For example, the tree associated with $243411$ is

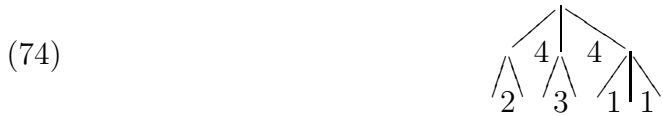
\begin{equation}
\vcenter{\xymatrix@C=0.5mm@R=4mm{
*{} & *{} & *{} & *{} & {}\ar@{-}[dlll]\ar@{-}[d]\ar@{-}[drrrr] \\
*{} &  {}\ar@{-}[dl]\ar@{-}[dr] & *{}
& *{4} & {}\ar@{-}[dl]\ar@{-}[dr] & *{}
& *{4} & *{} & {}\ar@{-}[dll]\ar@{-}[d]\ar@{-}[drr] \\ 
 {}
& *{2} &  {} &  {}
& *{3} &  {} &  {}
& *{1} &  {}
& *{1} &  {} \\
}}
\end{equation}

From the previous considerations, one can now
deduce a closed formula for the number of packed
words yielding a given plane tree, which can be regarded
as another generalization of the hook length
formula for binary trees:

\begin{theorem}
\label{Ft}
If a term $F_T(1)$ in the plane tree solution has the decomposition
\begin{equation}
F_T(1)=\sum_k c_k{t\choose k}
\end{equation}
then, $c_k$ is the number of packed words $u$ with maximal letter
$k$ such that $\TT(u)=T$.
\end{theorem}

\Proof A straightforward induction, from (\ref{wFk}) and (\ref{wT}).
\qed

For example, the following tree
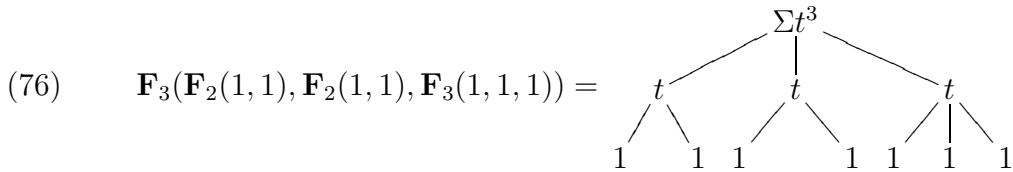
\entrymodifiers={+<4pt>}
\begin{equation}
\F_3(\F_2(1,1),\F_2(1,1),\F_3(1,1,1))=
\vcenter{\xymatrix@C=2mm@R=5mm{
*{} & *{} & *{} & *{} & {\Sigma t^3}\ar@{-}[dlll]\ar@{-}[d]\ar@{-}[drrrr] \\
*{} &  {t}\ar@{-}[dl]\ar@{-}[dr] & *{} & *{} 
    & {t}\ar@{-}[dl]\ar@{-}[dr] & *{} & *{} & *{}
    & {t}\ar@{-}[dll]\ar@{-}[d]\ar@{-}[drr] \\ 
 {1} & *{} &  {1} &  {1} & *{} &  {1} &  {1} & *{} &  {1} & *{} &  {1} \\
}}
\end{equation}
gives
\begin{equation}
\Sigma_0^t s^3\delta s=6{t\choose 4}+6{t\choose 3}+{t\choose 2}\,
\end{equation}
so that there are $6+6+1=13$ packed words whose plane trees have this shape:

%\entrymodifiers={+<4pt>[o][F]}
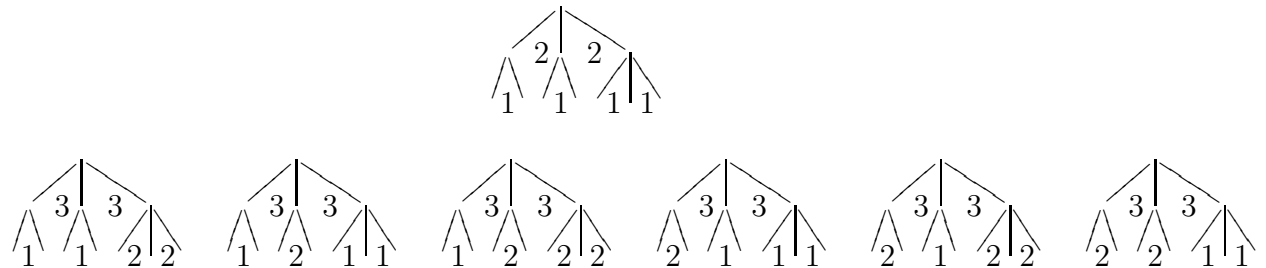
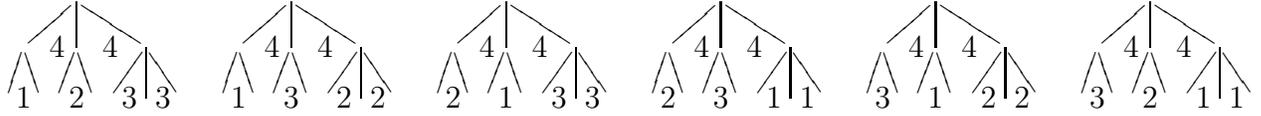
\entrymodifiers={+<4pt>}
\begin{equation*}
\vcenter{\xymatrix@C=0.5mm@R=4mm{
*{} & *{} & *{} & *{} & {}\ar@{-}[dlll]\ar@{-}[d]\ar@{-}[drrrr] \\
*{} &  {}\ar@{-}[dl]\ar@{-}[dr] & *{}
& *{2} & {}\ar@{-}[dl]\ar@{-}[dr] & *{}
& *{2} & *{} & {}\ar@{-}[dll]\ar@{-}[d]\ar@{-}[drr] \\ 
 {}
& *{1} &  {} &  {}
& *{1} &  {} &  {}
& *{1} &  {}
& *{1} &  {} \\
}}
\end{equation*}

\begin{equation*}
\vcenter{\xymatrix@C=0.5mm@R=4mm{
*{} & *{} & *{} & *{} & {}\ar@{-}[dlll]\ar@{-}[d]\ar@{-}[drrrr] \\
*{} &  {}\ar@{-}[dl]\ar@{-}[dr] & *{}
& *{3} & {}\ar@{-}[dl]\ar@{-}[dr] & *{}
& *{3} & *{} & {}\ar@{-}[dll]\ar@{-}[d]\ar@{-}[drr] \\ 
 {}
& *{1} &  {} &  {}
& *{1} &  {} &  {}
& *{2} &  {}
& *{2} &  {} \\
}}
\quad
\vcenter{\xymatrix@C=0.5mm@R=4mm{
*{} & *{} & *{} & *{} & {}\ar@{-}[dlll]\ar@{-}[d]\ar@{-}[drrrr] \\
*{} &  {}\ar@{-}[dl]\ar@{-}[dr] & *{}
& *{3} & {}\ar@{-}[dl]\ar@{-}[dr] & *{}
& *{3} & *{} & {}\ar@{-}[dll]\ar@{-}[d]\ar@{-}[drr] \\ 
 {}
& *{1} &  {} &  {}
& *{2} &  {} &  {}
& *{1} &  {}
& *{1} &  {} \\
}}
\quad
\vcenter{\xymatrix@C=0.5mm@R=4mm{
*{} & *{} & *{} & *{} & {}\ar@{-}[dlll]\ar@{-}[d]\ar@{-}[drrrr] \\
*{} &  {}\ar@{-}[dl]\ar@{-}[dr] & *{}
& *{3} & {}\ar@{-}[dl]\ar@{-}[dr] & *{}
& *{3} & *{} & {}\ar@{-}[dll]\ar@{-}[d]\ar@{-}[drr] \\ 
 {}
& *{1} &  {} &  {}
& *{2} &  {} &  {}
& *{2} &  {}
& *{2} &  {} \\
}}
\quad
\vcenter{\xymatrix@C=0.5mm@R=4mm{
*{} & *{} & *{} & *{} & {}\ar@{-}[dlll]\ar@{-}[d]\ar@{-}[drrrr] \\
*{} &  {}\ar@{-}[dl]\ar@{-}[dr] & *{}
& *{3} & {}\ar@{-}[dl]\ar@{-}[dr] & *{}
& *{3} & *{} & {}\ar@{-}[dll]\ar@{-}[d]\ar@{-}[drr] \\ 
 {}
& *{2} &  {} &  {}
& *{1} &  {} &  {}
& *{1} &  {}
& *{1} &  {} \\
}}
\quad
\vcenter{\xymatrix@C=0.5mm@R=4mm{
*{} & *{} & *{} & *{} & {}\ar@{-}[dlll]\ar@{-}[d]\ar@{-}[drrrr] \\
*{} &  {}\ar@{-}[dl]\ar@{-}[dr] & *{}
& *{3} & {}\ar@{-}[dl]\ar@{-}[dr] & *{}
& *{3} & *{} & {}\ar@{-}[dll]\ar@{-}[d]\ar@{-}[drr] \\ 
 {}
& *{2} &  {} &  {}
& *{1} &  {} &  {}
& *{2} &  {}
& *{2} &  {} \\
}}
\quad
\vcenter{\xymatrix@C=0.5mm@R=4mm{
*{} & *{} & *{} & *{} & {}\ar@{-}[dlll]\ar@{-}[d]\ar@{-}[drrrr] \\
*{} &  {}\ar@{-}[dl]\ar@{-}[dr] & *{}
& *{3} & {}\ar@{-}[dl]\ar@{-}[dr] & *{}
& *{3} & *{} & {}\ar@{-}[dll]\ar@{-}[d]\ar@{-}[drr] \\ 
 {}
& *{2} &  {} &  {}
& *{2} &  {} &  {}
& *{1} &  {}
& *{1} &  {} \\
}}
\end{equation*}

\begin{equation*}
\vcenter{\xymatrix@C=0.5mm@R=4mm{
*{} & *{} & *{} & *{} & {}\ar@{-}[dlll]\ar@{-}[d]\ar@{-}[drrrr] \\
*{} &  {}\ar@{-}[dl]\ar@{-}[dr] & *{}
& *{4} & {}\ar@{-}[dl]\ar@{-}[dr] & *{}
& *{4} & *{} & {}\ar@{-}[dll]\ar@{-}[d]\ar@{-}[drr] \\ 
 {}
& *{1} &  {} &  {}
& *{2} &  {} &  {}
& *{3} &  {}
& *{3} &  {} \\
}}
\quad
\vcenter{\xymatrix@C=0.5mm@R=4mm{
*{} & *{} & *{} & *{} & {}\ar@{-}[dlll]\ar@{-}[d]\ar@{-}[drrrr] \\
*{} &  {}\ar@{-}[dl]\ar@{-}[dr] & *{}
& *{4} & {}\ar@{-}[dl]\ar@{-}[dr] & *{}
& *{4} & *{} & {}\ar@{-}[dll]\ar@{-}[d]\ar@{-}[drr] \\ 
 {}
& *{1} &  {} &  {}
& *{3} &  {} &  {}
& *{2} &  {}
& *{2} &  {} \\
}}
\quad
\vcenter{\xymatrix@C=0.5mm@R=4mm{
*{} & *{} & *{} & *{} & {}\ar@{-}[dlll]\ar@{-}[d]\ar@{-}[drrrr] \\
*{} &  {}\ar@{-}[dl]\ar@{-}[dr] & *{}
& *{4} & {}\ar@{-}[dl]\ar@{-}[dr] & *{}
& *{4} & *{} & {}\ar@{-}[dll]\ar@{-}[d]\ar@{-}[drr] \\ 
 {}
& *{2} &  {} &  {}
& *{1} &  {} &  {}
& *{3} &  {}
& *{3} &  {} \\
}}
\quad
\vcenter{\xymatrix@C=0.5mm@R=4mm{
*{} & *{} & *{} & *{} & {}\ar@{-}[dlll]\ar@{-}[d]\ar@{-}[drrrr] \\
*{} &  {}\ar@{-}[dl]\ar@{-}[dr] & *{}
& *{4} & {}\ar@{-}[dl]\ar@{-}[dr] & *{}
& *{4} & *{} & {}\ar@{-}[dll]\ar@{-}[d]\ar@{-}[drr] \\ 
 {}
& *{2} &  {} &  {}
& *{3} &  {} &  {}
& *{1} &  {}
& *{1} &  {} \\
}}
\quad
\vcenter{\xymatrix@C=0.5mm@R=4mm{
*{} & *{} & *{} & *{} & {}\ar@{-}[dlll]\ar@{-}[d]\ar@{-}[drrrr] \\
*{} &  {}\ar@{-}[dl]\ar@{-}[dr] & *{}
& *{4} & {}\ar@{-}[dl]\ar@{-}[dr] & *{}
& *{4} & *{} & {}\ar@{-}[dll]\ar@{-}[d]\ar@{-}[drr] \\ 
 {}
& *{3} &  {} &  {}
& *{1} &  {} &  {}
& *{2} &  {}
& *{2} &  {} \\
}}
\quad
\vcenter{\xymatrix@C=0.5mm@R=4mm{
*{} & *{} & *{} & *{} & {}\ar@{-}[dlll]\ar@{-}[d]\ar@{-}[drrrr] \\
*{} &  {}\ar@{-}[dl]\ar@{-}[dr] & *{}
& *{4} & {}\ar@{-}[dl]\ar@{-}[dr] & *{}
& *{4} & *{} & {}\ar@{-}[dll]\ar@{-}[d]\ar@{-}[drr] \\ 
 {}
& *{3} &  {} &  {}
& *{2} &  {} &  {}
& *{1} &  {}
& *{1} &  {} \\
}}
\end{equation*}

\bigskip

\section{Functional equations associated to some generalizations of the
hook length formula}

\subsection{Postnikov's identity and  Eisenstein's exponential series}

Postnikov \cite{apost} has obtained the following identity
\begin{equation}\label{post}
(n+1)^{n-1} = \frac{n!}{2^n}\sum_{T\in\BT_n}
\prod_{v\in T}\left(1+\frac{1}{h_v}\right)\,.
\end{equation}
where $\BT_n$ is the set of (incomplete) binary trees with $n$
nodes.
Combinatorial proofs are given in \cite{chenyang,Seo}, and generalization
(to be discussed below) occur in \cite{DuLiu}.

Let $g(t)$ be the exponential generating function of the l.h.s
of (\ref{post}), that is,
\begin{equation}
g(t)=\sum_{n\ge 0}(n+1)^{n-1}\frac{t^n}{n!}\,.
\end{equation}
This is a famous power series, known as Eisenstein's 
generalized exponential (see, e.g., \cite{NTlag}). 
It satisfies the functional
equation
\begin{equation}\label{eis}
g(t)=e^{tg(t)}\,.
\end{equation}
Hence, $x=g(t)$ is solution of the differential equation
\begin{equation}
x' = x^2+txx' = x^2+t\frac{d}{dt}\left(\frac{x^2}{2}\right)\,,
\end{equation}
and integrating by parts, we obtain the fixed point equation
\begin{equation}
x = 1 + t\frac{x^2}{2}+\frac12\int_0^t x^2(s)ds = 1+B(x,x)
\end{equation}
with
\begin{equation}
B(x,y) =t\frac{xy}{2}+\frac12\int_0^t x(s)y(s)ds\,.
\end{equation}
From this, one derives that
%It is immediate to verify that in the binary tree expansion of the solution,
\begin{equation}
B_T(1)=\frac{1}{2^n}\prod_{v\in T}\left(1+\frac{1}{h_v}\right)t^n\,,
\end{equation}
since, by induction, if $T$ has $T_1$ (resp. $T_2$) as left (resp. right)
subtree with $n_1$ nodes (resp. $n_2$ nodes), then
\begin{equation}
\begin{split}
B_T(1) &= B(B_{T_1}(1), B_{T_2}(1))\\
& =
\frac{1}{2^{n_1}}\prod_{v\in T_1}\left(1+\frac{1}{h_v}\right)
\frac{1}{2^{n_2}}\prod_{v\in T_2}\left(1+\frac{1}{h_v}\right)
\left(\frac{1}{2}t^{n_1+n_2+1} + \frac{1}{2} \frac{t^{n_1+n_2+1}}{n_1+n_2+1}
\right) \\
& = 
\frac{1}{2^{n_1+n_2+1}}\prod_{v\in T}\left(1+\frac{1}{h_v}\right) t^n\,,
\end{split}
\end{equation}
which explains (\ref{post}). Note in particular that both terms of
$B(x,y)$ contribute to one term (either $1$ or $1/h_v$) for each node.

%%%%%%%%%%%%%%%%%%%%%%%%%%%%%%%%%%%%%%%%%%%%%%%%%%%%%%%%%%%%%%%%%%%%%%%%%%%%%%%
\subsection{Du-Liu identities}
Lascoux proposed a one parameter-generalization of (\ref{post}):
\begin{equation}\label{las1}
\sum_{T}\prod_{v}\left(\alpha +\frac{1}{h_v}\right)=\frac{1}{(n+1)!}\prod_{i
= 0}^{n-1} \left((n + 1 + i)\alpha  + n +1 -i\right).
\end{equation}
which has been proved by Du and Liu \cite{DuLiu}, who reformulated it as
\begin{equation}\label{las2}
\sum_{T}\prod_{v}\frac{(h_v+1)\alpha  + 1 -
h_v}{2h_v}=\frac{1}{n+1}{{(n+1)\alpha } \choose n}.
\end{equation}
and obtained the further generalization
\begin{equation}\label{las3}
\sum_{T}\prod_{v}\frac{(mh_v+1)\alpha +1-h_{v}}{(m+1)h_v}=\frac{1}{mn+1}{{(mn+1)\alpha }
\choose n}.
\end{equation}
where now, $T$ runs over plane $(m+1)$-ary trees.

These identities can also be  obtained from the tree solution of a functional
equation. Let $x=f(t)$ be the ordinary generating function of the
r.h.s. of (\ref{las3}), that is,
\begin{equation}
f(t)=\sum_{n\ge 0}{{(mn+1)\alpha } \choose n}\frac{t^n}{mn+1}\,.
\end{equation}
It follows from the Lagrange inversion formula (see, e.g., 
\cite[p. 35 ex.  25]{Mcd}) that $x$ is  solution of the
fixed point equation
\begin{equation}
x = (1+tx^m)^\alpha\,.
\end{equation}
Taking derivatives, we obtain the differential equation
\begin{equation}
x'=\alpha x^{m+1}+(\alpha m -1)t\frac{d}{dt}\left(\frac{x^{m+1}}{m+1}\right)
\end{equation}
and integrating by parts, we arrive at
\begin{equation}
x=1+\frac{\alpha m -1}{m+1} t x^{m+1}+\frac{\alpha+1}{m+1}
\int_0^tx^{m+1}(s)ds \equiv 1+F_{m+1}(x,x,\ldots,x)\,.
\end{equation}
As in the Postnikov identity, the $(m+1)$-ary tree expansion of the
solution associates to each tree $T$ the l.h.s. of (\ref{las3}), where both
terms of $F_{m+1}$ contribute to one term (either with coefficient $1$ or
$1/h_v$) for each node.

%%%%%%%%%%%%%%%%%%%%%%%%%%%%%%%%%%%%%%%%%%%%%%%%%%%%%%%%%%%%%%%%%%%%%%%%%%%%%%%
\section{Concluding remarks}

The original hook length formula for Young tableaux can be interpreted as
giving the image of a Schur function by the ring homomorphism
$f\mapsto f(\EE)$ defined on the power sums
\begin{equation}
 p_n \mapsto p_n(\EE) = \left\{
  \begin{array}{ll}
    1 & n=1 \\
    0 & n>1\\
  \end{array}\right.
\end{equation}
These are generalizations giving the images by the morphisms
\begin{equation}
\left\{
\begin{array}{l}
\displaystyle p_n \mapsto p_n(\frac{1}{1-q}) = \frac{1}{1-q^n},\\[11pt]
\displaystyle p_n \mapsto p_n(\alpha) = \alpha,\\[9pt]
\displaystyle p_n \mapsto p_n(\frac{1-t}{1-q}) = \frac{1-t^n}{1-q^n},\\
\end{array}\right.
%\begin{split}
%& p_n \mapsto p_n(\frac{1}{1-q}) = \frac{1}{1-q^n}, \qquad
% p_n \mapsto p_n(\alpha), \text{\ \ and} \\
%& p_n \mapsto p_n(\frac{1-t}{1-q}) = \frac{1-t^n}{1-q^n},
%\end{split}
\end{equation}
the last one giving back the first one for $t=0$ and the second one for
$t=q^\alpha$ and $q\to1$.

The theory of noncommutative symmetric functions allows one to define analogs
of these specializations for quasi-symmetric functions~\cite{NCSF2}, and
therefore also for those combinatorial Hopf algebras $H$ which admit
homomorphisms $H\to\QSym$.
This is the case of $\PBT$ and $\WQSym$, and Corollary~\ref{qhook} and
Theorem~\ref{Ft} can be interpreted as evaluation of $\P_T(1/(1-q))$ and
$\M_T(\alpha)$ respectively. It will be shown in a forthcoming paper that it
is in fact possible to evaluate both $\P_T$ and $\M_T$ on $(1-t)/(1-q)$
defined in the right way, and to get ($q$,$t$)-hook length formulas for binary
and plane trees.

%%%%%%%%%%%%%%%%%%%%%%%%%%%%%%%%%%%%%%%%%%%%%%%%%%%%%%%%%%%%%%%%%%%%%%%%%%%%%%%
%%%%%%%%%%%%%%%%%%%%%%%%%%%%%%%%%%%%%%%%%%%%%%%%%%%%%%%%%%%%%%%%%%%%%%%%%%%%%%%
%%%%%%%%%%%%%%%%%%%%%%%%%%%%%%%%%%%%%%%%%%%%%%%%%%%%%%%%%%%%%%%%%%%%%%%%%%%%%%%

\end{document}